\documentclass[12pt]{amsart}

\usepackage{amsfonts}
\usepackage{amssymb}
\usepackage{latexsym}
\usepackage{graphics}
\usepackage[2emode]{psfrag}
\usepackage{mathrsfs}
\usepackage{amsthm}
\usepackage{amsmath}
\usepackage[all]{xy}
\usepackage{enumerate}
\usepackage{hyperref}

\begin{document}

\newcommand{\thlabel}[1]{\label{th:#1}}
\newcommand{\thref}[1]{Theorem~\ref{th:#1}}
\newcommand{\selabel}[1]{\label{se:#1}}
\newcommand{\seref}[1]{Section~\ref{se:#1}}
\newcommand{\lelabel}[1]{\label{le:#1}}
\newcommand{\leref}[1]{Lemma~\ref{le:#1}}
\newcommand{\prlabel}[1]{\label{pr:#1}}
\newcommand{\prref}[1]{Proposition~\ref{pr:#1}}
\newcommand{\colabel}[1]{\label{co:#1}}
\newcommand{\coref}[1]{Corollary~\ref{co:#1}}
\newcommand{\relabel}[1]{\label{re:#1}}
\newcommand{\reref}[1]{Remark~\ref{re:#1}}
\newcommand{\exlabel}[1]{\label{ex:#1}}
\newcommand{\exref}[1]{Example~\ref{ex:#1}}
\newcommand{\delabel}[1]{\label{de:#1}}
\newcommand{\deref}[1]{Definition~\ref{de:#1}}
\newcommand{\eqlabel}[1]{\label{eq:#1}}
\newcommand{\equref}[1]{(\ref{eq:#1})}

\newcommand{\norm}[1]{\| #1 \|}
\def\*c{{}^*\hspace*{-1pt}{\cc}}
\def\C{\mathcal C}
\def\D{\mathcal D}
\def\J{\mathcal J}
\def\M{\mathcal M}
\def\T{\mathcal T}

\newcommand{\Ext}{{\rm Ext}}
\newcommand{\Fun}{{\rm Fun}}
\newcommand{\Mor}{{\rm Mor}\,}
\newcommand{\Aut}{{\rm Aut}\,}
\newcommand{\Hopf}{{\rm Hopf}\,}
\newcommand{\Ann}{{\rm Ann}\,}
\newcommand{\Ker}{{\rm Ker}\,}
\newcommand{\Coker}{{\rm Coker}\,}
\newcommand{\im}{{\rm Im}\,}
\newcommand{\coim}{{\rm Coim}\,}
\newcommand{\colim}{{\rm colim}\,}
\newcommand{\Trace}{{\rm Trace}\,}
\newcommand{\Char}{{\rm Char}\,}
\newcommand{\Spec}{{\rm Spec}\,}
\newcommand{\Span}{{\rm Span}\,}
\newcommand{\sgn}{{\rm sgn}\,}
\newcommand{\Id}{{\rm Id}\,}
\newcommand{\Com}{{\rm Com}\,}
\newcommand{\codim}{{\rm codim}}
\newcommand{\Mat}{{\rm Mat}}
\newcommand{\Coint}{{\rm Coint}}
\newcommand{\Incoint}{{\rm Incoint}}
\newcommand{\can}{{\sf can}}
\newcommand{\sign}{{\rm sign}}
\newcommand{\kar}{{\rm kar}}
\newcommand{\rad}{{\rm rad}}
\newcommand{\Rat}{{\rm Rat}}
\newcommand{\sd}{{\sf d}}

\def\Ab{\underline{\underline{\rm Ab}}}
\def\lan{\langle}
\def\ran{\rangle}
\def\ot{\otimes}
\def\cot{\square}
\def\bul{\bullet}
\def\ubul{\underline{\bullet}}
\def\tildej{\tilde{\jmath}}
\def\barj{\bar{\jmath}}

\def\id{\textrm{{\small 1}\normalsize\!\!1}}
\def\To{{\multimap\!\to}}
\def\bigperp{{\LARGE\textrm{$\perp$}}}
\newcommand{\QED}{\hspace{\stretch{1}}
\makebox[0mm][r]{$\Box$}\\}

\def\AA{{\mathbb A}}
\def\BB{{\mathbb B}}
\def\CC{{\mathbb C}}
\def\DD{{\mathbb D}}
\def\EE{{\mathbb E}}
\def\FF{{\mathbb F}}
\def\GG{{\mathbb G}}
\def\HH{{\mathbb H}}
\def\II{{\mathbb I}}
\def\JJ{{\mathbb J}}
\def\KK{{\mathbb K}}
\def\LL{{\mathbb L}}
\def\MM{{\mathbb M}}
\def\NN{{\mathbb N}}
\def\OO{{\mathbb O}}
\def\PP{{\mathbb P}}
\def\QQ{{\mathbb Q}}
\def\RR{{\mathbb R}}
\def\SS{{\mathbb S}}
\def\TT{{\mathbb T}}
\def\UU{{\mathbb U}}
\def\VV{{\mathbb V}}
\def\WW{{\mathbb W}}
\def\XX{{\mathbb X}}
\def\YY{{\mathbb Y}}
\def\ZZ{{\mathbb Z}}

\def\L{{\bf L}}
\def\G{{\bf G}}
\def\H{{\bf H}}
\def\F{{\bf F}}
\def\R{{\bf R}}
\def\P{{\bf P}}
\def\U{{\bf U}}

\def\aa{{\mathfrak A}}
\def\bb{{\mathfrak B}}
\def\cc{{\mathfrak C}}
\def\dd{{\mathfrak D}}
\def\ee{{\mathfrak E}}
\def\ff{{\mathfrak F}}
\def\gg{{\mathfrak G}}
\def\hh{{\mathfrak H}}
\def\ii{{\mathfrak I}}
\def\jj{{\mathfrak J}}
\def\kk{{\mathfrak K}}
\def\ll{{\mathfrak L}}
\def\mm{{\mathfrak M}}
\def\nn{{\mathfrak N}}
\def\oo{{\mathfrak O}}
\def\pp{{\mathfrak P}}
\def\qq{{\mathfrak Q}}
\def\rr{{\mathfrak R}}
\def\ss{{\mathfrak S}}
\def\tt{{\mathfrak T}}
\def\uu{{\mathfrak U}}
\def\vv{{\mathfrak V}}
\def\ww{{\mathfrak W}}
\def\xx{{\mathfrak X}}
\def\yy{{\mathfrak Y}}
\def\zz{{\mathfrak Z}}

\def\aaa{{\mathfrak a}}
\def\bbb{{\mathfrak b}}
\def\ccc{{\mathfrak c}}
\def\ddd{{\mathfrak d}}
\def\eee{{\mathfrak e}}
\def\fff{{\mathfrak f}}
\def\ggg{{\mathfrak g}}
\def\hhh{{\mathfrak h}}
\def\iii{{\mathfrak i}}
\def\jjj{{\mathfrak j}}
\def\kkk{{\mathfrak k}}
\def\lll{{\mathfrak l}}
\def\mmm{{\mathfrak m}}
\def\nnn{{\mathfrak n}}
\def\ooo{{\mathfrak o}}
\def\ppp{{\mathfrak p}}
\def\qqq{{\mathfrak q}}
\def\rrr{{\mathfrak r}}
\def\sss{{\mathfrak s}}
\def\ttt{{\mathfrak t}}
\def\uuu{{\mathfrak u}}
\def\vvv{{\mathfrak v}}
\def\www{{\mathfrak w}}
\def\xxx{{\mathfrak x}}
\def\yyy{{\mathfrak y}}
\def\zzz{{\mathfrak z}}

\newcommand{\aA}{\mathscr{A}}
\newcommand{\bB}{\mathscr{B}}
\newcommand{\cC}{\mathscr{C}}
\newcommand{\dD}{\mathscr{D}}
\newcommand{\eE}{\mathscr{E}}
\newcommand{\fF}{\mathscr{F}}
\newcommand{\gG}{\mathscr{G}}
\newcommand{\hH}{\mathscr{H}}
\newcommand{\iI}{\mathscr{I}}
\newcommand{\jJ}{\mathscr{J}}
\newcommand{\kK}{\mathscr{K}}
\newcommand{\lL}{\mathscr{L}}
\newcommand{\mM}{\mathscr{M}}
\newcommand{\nN}{\mathscr{N}}
\newcommand{\oO}{\mathscr{O}}
\newcommand{\pP}{\mathscr{P}}
\newcommand{\qQ}{\mathscr{Q}}
\newcommand{\rR}{\mathscr{R}}
\newcommand{\sS}{\mathscr{S}}
\newcommand{\tT}{\mathscr{T}}
\newcommand{\uU}{\mathscr{U}}
\newcommand{\vV}{\mathscr{V}}
\newcommand{\wW}{\mathscr{W}}
\newcommand{\xX}{\mathscr{X}}
\newcommand{\yY}{\mathscr{Y}}
\newcommand{\zZ}{\mathscr{Z}}

\newcommand{\Aa}{\mathcal{A}}
\newcommand{\Bb}{\mathcal{B}}
\newcommand{\Cc}{\mathcal{C}}
\newcommand{\Dd}{\mathcal{D}}
\newcommand{\Ee}{\mathcal{E}}
\newcommand{\Ff}{\mathcal{F}}
\newcommand{\Gg}{\mathcal{G}}
\newcommand{\Hh}{\mathcal{H}}
\newcommand{\Ii}{\mathcal{I}}
\newcommand{\Jj}{\mathcal{J}}
\newcommand{\Kk}{\mathcal{K}}
\newcommand{\Ll}{\mathcal{L}}
\newcommand{\Mm}{\mathcal{M}}
\newcommand{\Nn}{\mathcal{N}}
\newcommand{\Oo}{\mathcal{O}}
\newcommand{\Pp}{\mathcal{P}}
\newcommand{\Qq}{\mathcal{Q}}
\newcommand{\Rr}{\mathcal{R}}
\newcommand{\Ss}{\mathcal{S}}
\newcommand{\Tt}{\mathcal{T}}
\newcommand{\Uu}{\mathcal{U}}
\newcommand{\Vv}{\mathcal{V}}
\newcommand{\Ww}{\mathcal{W}}
\newcommand{\Xx}{\mathcal{X}}
\newcommand{\Yy}{\mathcal{Y}}
\newcommand{\Zz}{\mathcal{Z}}

\def\units{{\mathbb G}_m}
\def\rightact{\hbox{$\leftharpoonup$}}
\def\leftact{\hbox{$\rightharpoonup$}}

\def\text#1{{\rm {\rm #1}}}

\def\smashco{\mathrel>\joinrel\mathrel\triangleleft}
\def\cosmash{\mathrel\triangleright\joinrel\mathrel<}

\def\ol{\overline}
\def\ul{\underline}
\def\dul#1{\underline{\underline{#1}}}
\def\Nat{\dul{\rm Nat}}
\def\Set{\dul{\rm Set}}

\renewcommand{\subjclassname}{\textup{2000} Mathematics Subject
     Classification}

\swapnumbers
\newtheorem{proposition}{Proposition}[section]
\newtheorem{lemma}[proposition]{Lemma}
\newtheorem{corollary}[proposition]{Corollary}
\newtheorem{theorem}[proposition]{Theorem}

\theoremstyle{definition}
\newtheorem{definition}[proposition]{Definition}
\newtheorem{example}[proposition]{Example}
\newtheorem{examples}[proposition]{Examples}

\theoremstyle{remark}
\newtheorem{remarks}[proposition]{Remarks}
\newtheorem{remark}[proposition]{Remark}


\newcommand{\cat}[1]{\mathcal{#1}}
\newcommand{\unidad}[1]{\eta_{#1}}
\newcommand{\counidad}[1]{\delta_{#1}}
\newcommand{\stat}[1]{\mathrm{Stat}(#1)}
\newcommand{\adstat}[1]{\mathrm{Adst}(#1)}
\newcommand{\coker}[1]{\mathrm{coker}(#1)}
\newcommand{\rcomod}[1]{\mathcal{M}^{#1}}
\newcommand{\lcomod}[1]{{}^{#1}\mathcal{M}}
\newcommand{\rmod}[1]{\mathcal{M}_{#1}}
\newcommand{\lmod}[1]{{}_{#1}\mathcal{M}}
\newcommand{\coend}[2]{\mathrm{Coend}_{#1}(#2)}
\newcommand{\End}[3]{\mathrm{End}_{#1}^{#2}(#3)}
\renewcommand{\hom}[3]{\mathrm{Hom}_{#1}(#2,#3)}
\newcommand{\Hom}{\mathrm{Hom}}
\newcommand{\cotensor}[1]{\square_{#1}}
\newcommand{\functor}[1]{\mathbf{#1}}
\newcommand{\Mod}{\mbox{{\rm -Mod}}}
\newcommand{\cohom}[3]{\mathrm{h}_{#1}(#2,#3)}
\newcommand{\HOM}[3]{\mathrm{HOM}_{#1}(#2,#3)}

\title{Quasi-Frobenius functors. Applications}
\author[F.\ Casta\~ no Iglesias]{F.~Casta\~no Iglesias$^{\dagger}$}
\thanks{${^\dagger}$ Research partially
supported by Grant P07-FQM-03128 of Junta de Andalucia.}
\address{Departamento de
Estad\'{\i}stica y Matem\'atica Aplicada, Universidad
de Almer\'{\i}a 04120, Almer\'{\i}a, Spain}
\email{fci@ual.es}

 \author[C.\ N\v ast\v astescu]{C.~N\v{a}st\v{a}sescu${^\ddagger}$}
\thanks{${^\ddagger}$ Research
supported by Grant 434/1.10.2007 of CNCSIS, PN II (ID\_1005).}
\address{Facultatea de Matematic{\u a}, Str. Academiei 14 \\
 RO-70109, Bucharest,
 Romania}
\email{constantin\_nastasescu@yahoo.com}
\author[J.\ Vercruysse]{J.\ Vercruysse$^\diamond$}
\address{Faculty of Engineering, Vrije Universiteit Brussel (VUB), B-1050 Brussels, Belgium}
\email{jvercruy@vub.ac.be}
\urladdr{homepages.vub.ac.be/\~{}jvercruy}
\thanks{$^\diamond$ The author is Postdoctoral Fellow of the Fund for Scientific Research--Flanders
(Belgium) (F.W.O.--Vlaanderen).}
\date{}

\begin{abstract} We investigate functors between abelian categories having  a left adjoint and a right adjoint that are \emph{similar} (these functors are called \emph{quasi-Frobenius
functors}). We introduce the notion of a \emph{quasi-Frobenius bimodule} and give a characterization of these bimodules in terms of quasi-Frobenius functors. Some applications to corings and graded rings are presented. In particular, the concept of quasi-Frobenius homomorphism of corings is introduced. Finally, a version of the endomorphism ring Theorem for quasi-Frobenius extensions in terms of corings is obtained.
\end{abstract}

\maketitle


\section*{Introduction}
For any ring extension $\varphi: R\rightarrow S$, we can always consider the triple of functors $\Gamma = (S\otimes_R-, \varphi_*, \hom{R}{_RS}{-})$ where  $S\otimes_R-$ and $\hom{R}{_RS}{-}$ are, respectively, the left and the right adjoint of the restriction of scalars functor $\varphi_*: \lmod{S}\rightarrow \lmod{R}.$ The functor $\varphi_*$ is termed a \emph{Frobenius} functor if $S\otimes_R-$ and $\hom{R}{_RS_S}{-}$ are naturally isomorphic. Morita observed in \cite{Morita} that $\varphi_*$ is a Frobenius functor if and only if $\varphi$ is a Frobenius extension in the sense of \cite{Kas}, i.e.\ $S$ is finitely generated and projective as a left $R$-module and $S\cong \hom{R}{_RS}{R}$ as an $(S,R)$-bimodule. Frobenius algebras, and Frobenius ring extensions have been studied for more than hundred years, and surfaced in the most diverse parts of algebra and beyond.

The advantage of studying Frobenius functors, rather than Frobenius algebras, is twofold. First, it provides a more abstract, and therefore in many cases a more clarifying view on many aspects and properties of Frobenius algebras and ring extensions. Secondly, Frobenius functors can be studied as well between other than pure module categories. This made it possible to transfer properties of Frobenius ring extensions to very similar properties in terms of for example Hopf algebras and corings by studying Frobenius functors for their categories of Hopf modules or comodules, respectively.

Beside Frobenius algebras, there is also a vivid interest in quasi-Frobenius algebras. These are more general, but still posses many of the interesting properties of Frobenius algebras. The aim of this note is to initiate a functorial study of quasi-Frobenius algebras, similar to the Frobenius case.

Recall from M\"{u}ller \cite{muller} the notion of a left (and right) quasi-Frobenius extension, generalizing a Frobenius extension: $\varphi$ is a \emph{left quasi-Frobenius extension} if $_RS$ is finitely generated and projective and  $S$ is isomorphic as an $(S,R)$-bimodule to a direct summand of a finite direct sum of copies of $ \hom{R}{_RS}{R}$. Equivalently, $_RS$ and $S_R$ are finitely generated and projective and $\hom{R}{S_R}{R}$ is isomorphic as $(R,S)$-bimodule to a direct summand of a finite direct sum of copies of $S$. Similarly, $\varphi$ is a \emph{right quasi-Frobenius extension} if $S_R$ is finitely generated and projective and $S$ is isomorphic as $(R,S)$-bimodule to a direct summand of a finite direct sum of $\hom{R}{S_R}{R}$. We easily conclude that $\varphi$ is a \emph{quasi-Frobenius extension} (both left and right) if and only if $_RS$  is finitely generated and projective and the $(S,R)$-bimodules $S$ and $\hom{R}{_RS}{R}$ are similar, i.e.\ $(S\otimes_R-,\varphi_*)$ is a \emph{quasi-strongly adjoint pair} in the sense of \cite[Theorem 5.1]{Morita}. Note that $\varphi_*$ is a member of the triple of functors $\Gamma$ where  $S\otimes_R-$ and $\hom{R}{_RS_S}{-}$ are similar.  In this note, we will term a functor $\varphi_*$ with these properties a \emph{quasi-Frobenius functor}. So the ring extension $\varphi$ is quasi-Frobenius if and only if $\varphi_*$ is a quasi-Frobenius functor.

In this paper we shall concentrate on quasi-Frobenius functors between Grothendieck categories.
The purpose is to give a categorical framework to study quasi-Frobenius properties,
in a way that the results can be applied not only to ring extensions, but also to bimodules, graded rings, coring homomorphisms and bicomodules.
Moreover the functorial language can give a natural interpretation to certain known results, for example, M\"uller's result that if $\varphi:R\to S$ is a quasi-Frobenius ring extension of two $k$-algebras over a field, then $R$ is a quasi-Frobenius $k$-algebra if and only if $S$ is quasi-Frobenius $k$-algebra can now be understood as a consequence of our theory (see \coref{frobext}).

After recalling some elementary definitions, we introduce in \seref{functors} the concept of a \emph{quasi-Frobenius triple of functors} $(\L,\F,\R)$, or equivalently, a quasi-Frobenius functor $\F$, for categories with finite coproducts (see \deref{qffunctor}).
This concept generalizes the notion of a Frobenius functor \cite{CGN} and is closely related to the recently introduced left quasi-Frobenius pair of functors \cite{guo}. We prove some basic properties of quasi-Frobenius triples, in particular we show that they behave very well with respect to limits and colimits (see Lemma \ref{properties}).

In \seref{modules}, we characterize quasi-Frobenius functors between module categories, in fact, this coincides with Morita's notion of a strongly quasi-adjoint pair of functors. A particulary interesting feature of this special situation is that the left and right adjoint of a quasi-Frobenius functor is again quasi-Frobenius, hence it induces a sequence of quasi-Frobenius functors (see \reref{sequence}), a property that seems to be lost in the general case.

Another interesting case is given by graded rings and modules and this is considered in \seref{graded}, where we provide an example of a functor that is quasi-Frobenius if and only if it is Frobenius (see \reref{qf=f}).

The notion of Frobenius extension for coalgebras over fields was introduced by G\'omez-Torrecillas and the first two authors in \cite{CGN} and extended to corings by Zarouali-Darkaoui \cite{Za}. In \seref{comodules}, we characterize quasi-Frobenius functors between categories of comodules over corings. In the last section we introduce the concept of a \emph{quasi-Frobenius} morphism of corings which generalizes the notion of a Frobenius morphism of corings. Next we focus on corings for which the induction functor from the category of right $A$-modules to the category of right comodules of an $A$-coring is quasi-Frobenius. We term such corings \emph{quasi-Frobenius corings}. As an application, we prove that any quasi-Frobenius extension $\varphi:R\rightarrow S$ induces a quasi-Frobenius coring. This result generalizes \cite[Theorem 2.7]{Brze1} and \cite[Proposition 4.3]{guo}.

\section{Preliminaries}

\subsection*{Divisibility and similarity}
Let $\Aa$ be a category with finite coproducts. For an object $X$ in $\Aa$, and a positive integer $n$, we denote by $X^n$ the coproduct (direct sum) of $n$ copies of $X$. Consider now two objects $X$ and $Y$ in $\Aa$. Then we say that \emph{$X$ divides $Y$}, denoted by $X~|Y$, if there exist positive integer $n$ and morphisms $\phi:X\to Y^n$ and $\psi:Y^n\to X$, such that $\psi\circ\phi=1_X$. Two objects $X$ and $Y$ in $\Aa$ are said to be \emph{similar}, denoted by $X\sim Y$ if $X$ divides $Y$ and $Y$ divides $X$ at the same time.
Clearly, ``$\sim$'' defines an equivalence relation on the class of objects of $\Aa$.
From the splitting lemma, the following result follows immediately.
\begin{lemma}\lelabel{splitting}
Let $\Aa$ be an abelian category and $X, Y\in\Aa$ two objects. Then the following statements are equivalent
\begin{enumerate}[(i)]
\item $X~|~Y$;
\item there exists a positive integer $n$ and an object $Z\in \Aa$ such that $0\to Z\to Y^n\to X\to 0$ is a split exact sequence;
\item there exists a positive integer $n$ and an object $Z\in \Aa$ such that $Y^n\cong X\oplus Z$.
\end{enumerate}
\end{lemma}
In this paper we will consider some particular cases of this situation. In this section we introduce the necessary notation.

\subsection*{Functor categories}
Let $\Aa$ an $\Bb$ be two categories, where $\Aa$ has finite coproducts. Consider now a suitable category $\Fun(\Bb,\Aa)$, whose objects are functors $\F:\Bb\to \Aa$ and whose morphisms are natural transformations between these functors. One can easily check that $\Fun(\Bb,\Aa)$ has finite coproducts, which can be defined as follows. For any two functors $\F,\G:\Bb\to \Aa$, and $X\in\Bb$, we put $(\F\oplus \G)(X)=\F(X)\oplus \G(X)$. The definitions of divisibility and similarity therefore apply to $\Fun(\Aa,\Bb)$. Explicitly,
given covariant functors $\functor{L,R}:
\mathcal{B}\rightarrow \mathcal{A}$, we
say that $\functor{L}$ \emph{divides} $\functor{R}$, denoted by
$\functor{L}\mid\functor{R},$  if for some positive integer $n$
there are natural transformations
\[\xymatrix{\functor{L}(X)\ar[r]^-{\phi(X)}&  \functor{R}(X)^n
 \ar[r]^-{\psi(X)}& \functor{L}(X) },\]
such that $\psi(X)\circ\phi(X) = 1_{\functor{L}(X)}$, for all $X\in
\mathcal{B}$. Analogously, $\functor{R}\mid\functor{L}$ if for some
positive integer $m$ there are natural transformations
\[\xymatrix{\functor{R}(X)\ar[r]^-{\phi'(X)}&  \functor{L}(X)^m
 \ar[r]^-{\psi'(X)}& \functor{R}(X) },\]
such that $\psi'(X)\circ\phi'(X) = 1_{\functor{R}(X)}$, for all
$X\in \mathcal{B}$. The functor $\functor{L}$ is said to be
 \emph{similar} to $\functor{R}$, denoted by
 $\functor{L}\sim\functor{R},$ whenever  $\functor{L}\mid\functor{R}$ and
$\functor{R}\mid\functor{L}$.

The following lemma is obvious,
\begin{lemma}\lelabel{fundivobjdiv}
Consider functors $\L,\R:\Bb\to\Aa$. If $\L~|~\R$, then $\L(X)~|~\R(X)$ in $\Aa$ for all $X\in \Bb$. Consequently, if $\L\sim\R$, then $\L(X)\sim\R(X)$ in $\Aa$ for all $X\in\Bb$.
\end{lemma}

We will study properties of these functors in more detail in \seref{functors}.

\subsection*{Rings and modules}
Consider associative and unital rings  $R$ and $S$. We use the following notation and definitions as in \cite{A-F}.
For an $(R,S)$-bimodule we write $_RM$, $M_S$ or $_RM_S$ if we want to stress the left $R$-module, right $S$-module or bimodule structure on $M$, respectively. If $N$ is another $(R,S)$-bimodule, then $\hom{R}{M}{N}=
\hom{R}{_RM}{_RN}$ denotes the abelian group of left $R$-module maps, $\hom{S}{M}{N}$ is the abelian group of right $S$-module maps. A bimodule has two duals, a left dual $(_RM)^*=\Hom_R({_RM},{_RR})$ and a right dual $(M_S)^*=\Hom_S(M_S,S_S)$.
The categories of all $(R,S)$-bimodules, left $R$-modules and right $S$-modules are denoted respectively by $_R\mathcal{M}_S$, ${_R\Mm}$ and $\Mm_S$.

We can consider divisibility and similarity in ${_R\Mm_S}$. This leads now to the following explicit description.
We say that the $(R,S)$-bimodule $M$ devides the $(R,S)$-bimodule $N$, denoted by $M~|~N$, if there exists a positive integer $n$ and an $(R,S)$-bimodule $P$ such that $M\oplus P\cong N^{(n)}$.
Furthermore, $M$ and $N$ are called
\emph{similar}, abbreviated $_RM_S\sim \, _RN_S$, if $M$ devides $N$ and $N$ devides $M$, i.e.\ there are
$m,n\in \mathbb{N}$ and $(R,S)$-bimodules $P$ and $Q$ such that
$M\oplus P \cong N^{(n)}$ and $N\oplus Q \cong M^{(m)}$ as bimodules.
The properties of similar bimodules are examined in \seref{modules}.

\subsection*{Corings and comodules}

Let $A$ be an associative and unitary algebra over a commutative ring (with unit)
$k$. We recall from \cite{Sw2} that an $A$-coring $\mathfrak{C}$
consists of an $A$-bimodule $\mathfrak{C}$ with two $A$-bimodule
maps
$$\Delta: \mathfrak{C}\rightarrow \mathfrak{C}\otimes_A\mathfrak{C},
\mbox{      } \epsilon: \mathfrak{C}\rightarrow A$$ such that
$(\mathfrak{C}\otimes_A\Delta)\circ \Delta =
(\Delta\otimes_A\mathfrak{C})\circ \Delta$ and
$(\mathfrak{C}\otimes_A\epsilon)\circ \Delta =
(\epsilon\otimes_A\mathfrak{C})\circ \Delta = 1_{\mathfrak{C}}.$ A
right $\mathfrak{C}$-comodule is a pair $(M,\rho_M)$ consisting of a
right $A$-module $M$ and an $A$-linear map $\rho_M:M\rightarrow
M\otimes_A\mathfrak{C}$ satisfying $(M\otimes_A\Delta)\circ \rho_M =
(\rho_M\otimes_A\mathfrak{C})\circ \rho_M$ and
$(M\otimes_A\epsilon)\circ \rho_M =1_M$. Let $M$ and $N$ be two right $\cc$-comodules.
A morphism of right $\cc$-comodules, or a right $\cc$-colinear map, is a right $A$-linear map $f:M\to N$ such that $\rho_N\circ f= (f\ot_A\cc)\circ\rho_M$.
The right
$\mathfrak{C}$-comodules together with their morphisms form the
additive category $\mathcal{M}^{\mathfrak{C}}$. If $_A\mathfrak{C}$
is flat, then $\mathcal{M}^{\mathfrak{C}}$ is a Grothendieck
category.

Consider now an $A$-coring $\cc$ and a $B$-coring $\dd$. The category ${^\cc\Mm^\dd}$ consists of objects, called $(\cc,\dd)$-bicomodules, that are at the same time a left $\cc$-comodule $(M,\lambda_{M,\cc})$ and a right $\dd$-comodule $(M,\rho_{M,\dd})$ such that
$$(\cc\ot_A\rho_{M,\dd})\circ\lambda_{M,\cc}=(\lambda_{M,\cc}\ot_B\dd)\circ\rho_{M,\dd},$$
and whose morphisms are $k$-linear maps that are at the same time left $\cc$-colinear and right $\dd$-colinear.
Furthermore, given bicomodules $N\in
 \, ^{\mathfrak{D}}\rcomod{\mathfrak{C}}$ and $\overline{N}\in
 \, ^{\mathfrak{C}}\rcomod{\mathfrak{D}'}$ with $\mathfrak{D}'$ any
 $B'$-coring, we can consider the \emph{cotensor product} as the following equalizer in ${_B\Mm_{B'}}$.
\[
\xymatrix{
N\cotensor{\mathfrak{C}}\overline{N} \ar[rr] && M\ot_A \ol{N} \ar@<.5ex>[rr]^-{M\ot_A\lambda_{N,\cc}} \ar@<-.5ex>[rr]_-{\rho_{M,\cc}\ot_A\ol{N}} &&
M\ot_A\cc\ot_A\ol{N}.
}
\]
If $\mathfrak{D}_B$ and
 $_{B'}\mathfrak{D}$ are flat modules, then
 $N\cotensor{\mathfrak{C}}\overline{N}$ is a
 $(\mathfrak{D},\mathfrak{D}')$-bicomodule. In particular, $N\cotensor{\mathfrak{C}}\mathfrak{C}\cong
 N$ as $(\mathfrak{D},\mathfrak{C})$-bicomodule.

Applying the definitions of divisibility and similarity to the category ${^\dd\Mm^\cc}$, we obtain the following.
A bicomodule $_{\mathfrak{D}}N_{\mathfrak{C}}$ will be said to \emph{divide} the
bicomodule $_{\mathfrak{D}}\overline{N}_{\mathfrak{C}}$, abbreviated
$_{\mathfrak{D}}N_{\mathfrak{C}} ~|~ {_{\mathfrak{D}}\overline{N}_{\mathfrak{C}}}$, if there is an integer $n\in
\mathbb{N}$ and a $(\mathfrak{D},\mathfrak{C})$-bicomodule $P$ such that $N\oplus P \cong \overline{N}^{(n)}$.
Furthermore, $N$ is similar to $\ol{N}$, or $N\sim\ol{N}$ if in addition there is an integer $m\in\NN$ and a $(\dd,\cc)$-bicomodule $Q$ such that
$\overline{N}\oplus Q \cong N^{(n)}$ as bicomodules.

Similar bicomodules will be studied in \seref{comodules} and \seref{corings}.

\section{Quasi-Frobenius functors between Grothendieck categories}\selabel{functors}

Let $\Aa$ and $\Bb$ be categories that posses finite coproducts.
Consider a triple of functors $\Gamma
=(\functor{L,F,R})$, where $ \functor{F}:\mathcal{A}\rightarrow
\mathcal{B}$ has a left adjoint $\functor{L}: \mathcal{B}\rightarrow
\mathcal{A}$ and also a right adjoint
$\functor{R}:\mathcal{B}\rightarrow \mathcal{A}$. In this situation, we call $\Gamma$ an \emph{adjoint} triple. Notice that
$\functor{F}$ is exact and preserves limits and colimits,
$\functor{L}$ always preserves projective objects and colimits and is right exact, the functor $\functor{R}$ preserves injective objects and limits and is left exact.
\begin{definition}\delabel{qffunctor}
A \emph{quasi-Frobenius triple} for the categories $\mathcal{A}$ and $ \mathcal{B}$
consists of an adjoint triple of functors $\Gamma =(\functor{L,F,R})$ as above, where
$\functor{L}$ and $\functor{R}$ are similar functors.

A functor $ \functor{F}:\mathcal{A}\rightarrow \mathcal{B}$ is said
to be a \emph{quasi-Frobenius functor} if $(\functor{L,F,R})$ is a
quasi-Frobenius triple for some functors $\functor{L,R}:
\mathcal{B}\rightarrow \mathcal{A}$.
\end{definition}

\begin{remark}
In the case where $\Bb=\Mm_B$ and $\Aa=\Mm_A$ for certain rings $A$ and $B$, the situation of a quasi-Frobenius triple of functors for the categories $\Aa$ and $\Bb$ was termed a strongly adjoint pair of functors in \cite{Morita}. We will study this particular situation in \seref{modules}.
\end{remark}

Before proving some properties of quasi-Frobenius triples, we state the following elementary lemma.

\begin{lemma}\lelabel{limits}
Let $\Aa$ be a Grothendieck category. For a positive integer $n$, let $\P_n:\Aa\to \Aa$ be the defined by $\P_n(X)=X^n$. Then $\P_n$ preserves all limits and colimits.
\end{lemma}

\begin{proof}
Let us proof that $\P_n$ preserves arbitrary colimits, the proof for limits is obtained by applying dual arguments. Consider a (small) category $\Zz$ and a functor $\H:\Zz\to \Aa$. Let $(C,c_Z)=\colim \H$ be the colimit of $\H$. This means that $(C,c_Z)$ is a cocone on $\H$ (i.e. $C\in\Aa$ is an object and $c_Z:\H(Z)\to C$ is a collection of morphisms defined for all $Z\in\Zz$, such that for all $f:Z\to Z'$ in $\Zz$, we have $c_Z=c_{Z'}\circ \H(f)$), such that for every other cocone $(M,m_Z)$ there is a unique morphism $f:C\to M$ such that $m_Z=f\circ c_Z$ for all $Z\in\Zz$.

Now consider the functor $\P_n\circ \H:\Zz\to \Aa$, we have to show that $\P_n (\colim \H)=\colim (\P_n\circ \H)$. Let $(M,m_Z)$ be any cocone on $\P_n\circ \H$. Denote by $\pi_i:\H(Z)^n\to \H(Z)$ the projection on the $i$-th component in the direct sum. Then $(M,\pi_i\circ m_Z)$ is a cocone on $\H$ for any $i=1,\ldots,n$. Hence there exist unique morphisms $f_i:C\to M$ satisfying $\pi_i\circ m_Z=f_i\circ c_Z$ for all $i$ and all $Z$. In this way we obtain a unique morphism $f=f_1\oplus f_2 \oplus \ldots \oplus f_n : C^n\to M$ satisfying $m_Z=f\circ \P_n(c_Z)$. Therefore, $(\P_n(C^n),\P_n(c_Z))$ is the colimit of $\P_n\circ \H$.
\end{proof}

Remark that with notation as in \leref{limits}, $\F^n=\P_n\circ\F$ for any functor $\F:\Aa\to\Aa$.

\begin{lemma} \label{properties} Consider Grothendieck categories $\mathcal{A}, \mathcal{B}$
and $\mathcal{C}$. If $(\functor{L,F,R})$ is a quasi-Frobenius
triple for $\mathcal{A}$ and $ \mathcal{B}$, then
\begin{enumerate}[(a)]
\item The functors $\functor{L,F}$ and $\functor{R}$ are exact and preserve
all limits and colimits as well as injective and projective objects. The functors $\L$ and $\F$ preserve finitely generated objects.
\item If $( \functor{\bar L,
\bar F, \bar R})$ is  a  quasi-Frobenius  triple for $\mathcal{B}$ and $
\mathcal{C}$,
 then   $$(\functor{L\circ \bar L,\bar F\circ F,R\circ \bar R})$$ is also a
quasi-Frobenius triple for $\mathcal{A}$ and $\mathcal{C}$.
\end{enumerate}
\end{lemma}
\begin{proof}
$(a)$ Since $\functor{F}$ is exact, $\functor{L}$ preserves
projective objects and $\functor{R}$ injective objects. Suppose that
$\functor{L}\sim \functor{R}$. Then there are natural
transformations
\begin{equation}\label{ldivider}
\xymatrix{\functor{L}\ar[r]^-{\phi}&  \functor{R}^n
 \ar[r]^-{\psi}& \functor{L} }
 \end{equation}
such that $\psi\circ\phi = 1_{\functor{L}}$ and
\begin{equation}\label{rdividel}
\xymatrix{\functor{R}\ar[r]^-{\phi'}&  \functor{L}^m
 \ar[r]^-{\psi'}& \functor{R} }\end{equation}
such that $\psi'\circ\phi' = 1_{\functor{R}}$. From
$(\ref{ldivider})$ and $(\ref{rdividel})$, we easily deduce that
$\functor{R}$ preserves projective objects and  $\functor{L}$
 injective objects. Consider now any short exact sequence in
$\mathcal{B}$
\begin{equation}\label{exact}
\xymatrix{ 0\ar[r] & X\ar[r]^f &Y\ar[r]^g&Z\ar[r]&0} \end{equation}
 Applying $\functor{R}^n$ and $\L$
to (\ref{exact}), we obtain the commutative diagram with exact rows
\[
\xymatrix{0 \ar[r] & \functor{R}^n(X) \ar[r]^{\functor{R}^n(f)} & \functor{R}^n(Y) \ar[r]^{\functor{R}^n(g)} & \functor{R}^n(Z) &  \\
& \functor{L}(X) \ar[r]^{\functor{L}(f)} \ar[u]^{\phi_X} &
 \functor{L}(Y) \ar[r]^{\functor{L}(g)} \ar[u]^{\phi_Y}
 & \functor{L}(Z)  \ar[u]^{\phi_Z} \ar[r] &  0 \\ & 0 \ar[u] & 0 \ar[u] &0 \ar[u] &
 }
\]
from which it follows that $\functor{L}(f)$ is monic and therefore,
$\functor{L}$ exact. A similar argument with $\functor{L}^m$ shows
that  $\functor{R}$ is exact.

Since $\L$ is a left adjoint, we already know that it preserves colimits, similarly, $\R$, being a right adjoint, preserves limits. Before we show that $\R$ preserves arbitrary colimits, let us remark that by \leref{limits}, $\L^m$ and $\R^m$ preserve colimits and limits respectively. Now let $\Zz$ be a (small) category and $\H:\Zz\to \Bb$ a covariant functor.
We denote as in \leref{limits} $(C,c_Z)=\colim \H$. Obviously, $(\R(C),\R(c_Z))$ and $(\L^m(C),\L^m(c_Z))$ are cocones on respectively $\R\H$ and $\L^m\H$. Let us denote the colimit of $\R\H$ by $(C_R,r_Z)$ and the colimit of $\L^m\H$ by $(C_L,l_Z)$. Then there are unique morphisms $r:C_R\to \R(C)$ and $l:C_L\to \L^m(C)$ such that $\R(c_Z)=r\circ r_Z$ and $\L^m(c_Z)= l\circ l_Z$. That $\L^m$ preserves colimits means exactly that $l$ is an isomorphism, we have to show that $r$ is an isomorphism as well. From the properties of a colimit, one easily obtains that $r$ is a monomorphism.
Next, observe that $(C_R, r_Z\circ \psi'_{\H Z})$ is a cocone on $\L^m\H$ and $(C_L,l_Z\circ \phi'_{\H Z})$ is a cocone on $\R \H$. This induces morphisms $\bar\psi: C_L\to C_R$ and $\bar\phi:C_R\to C_L$ such that
$$r_Z\circ\psi'_{\H Z}=\bar\psi\circ l_Z \quad{\rm and}\quad
  l_Z\circ \phi'_{\H Z}=\bar\phi\circ r_Z.$$
Therefore, we find for all $Z\in\Zz$,
\begin{eqnarray*}
\phi'_C\circ r\circ r_Z&=& \phi'_C\circ \R(c_Z) = \L^m(c_Z)\circ \phi'_{\H Z}\\
&=& l\circ l_Z\circ \phi'_{\H Z} = l\circ \bar\phi\circ r_Z,
\end{eqnarray*}
where we used the naturality of $\phi'$ in the second equality. Since this holds for all $Z\in \Zz$, it follows from the properties of the colimit that $\phi'_C\circ r=l\circ \bar\phi$. A similar computation shows that $\psi'_C\circ l=r\circ \bar\psi$. Hence, we obtain the following commutative diagram where the upper row is exact,
\[
\xymatrix{
0\ar[rr] && \R(C) \ar[rr]^-{\phi'_C} && \L^m(C) \ar[rr]^-{\psi'_C} \ar[rr] && \R(C) \ar[rr] && 0 \\
&& C_R \ar[u]^{r} \ar[rr]^-{\bar\phi} && C_L \ar[u]^{l} \ar[rr]^-{\bar\psi} && C_R \ar[u]^{r}
}
\]
From this diagram one obtains that $r$ is an epimorphism, hence also an isomorphism.

Recall that an object $M$ of $\Aa$
is finitely generated if the functor $\hom{\Aa}{M}{-}$
preserves the sum of an arbitrary directed system $\{X_i\}_{i\in I}$ in $\Aa$.
By applying the adjunction property of $\L$ and $\F$, we find for a finitely generated object $M\in\Bb$ that
\begin{eqnarray*}
\Hom_\Aa(\L (M),\sum_{i\in I} X_i) &=& \Hom_\Bb(M,\F(\sum_{i\in I}X_i))=\Hom_\Bb(M,\sum_{i\in I}\F (X_i))\\
&=&\sum_{i\in I}\Hom_\Bb(M,\F (X_i))=\sum_{i\in I}\Hom_\Bb(\L (M), X_i),
\end{eqnarray*}
where we used the fact that $\F$ preserves colimits, hence sums, in the second equation and that $M$ is finitely generated in $\Bb$ in the third equation. It follows that $\L M$ is a finitely generated object in $\Aa$. Similarly, $\F N$ is finitely generated for a finitely generated object $N\in \Aa$.\\
$(b)$ Consider natural transformations $\phi$ and $\psi$ as in \eqref{ldivider} and \eqref{rdividel}, and similarly natural transformations $\bar\phi$ and $\bar\psi$ for the triple $(\bar\L,\bar\F,\bar\R)$.
Since by part $(a)$, all functors preserve in particular (finite) direct sums, we can define natural transformations
\begin{eqnarray*}
\xymatrix{\L \L' X \ar[r]^-{\L\bar{\phi}_X} & \L\bar\R X^{n'} \ar[r]^-{\phi_{\bar \R X^{n'}}} & \R\bar\R X^{n+n'}} ,
&&
\xymatrix{\R\bar\R X^{n+n'} \ar[r]^-{\psi_{\bar \R X^{n'}}} & \L\bar \R X^{n'} \ar[r]^-{\L\bar\psi_X} & \L\bar\L X};\\
\xymatrix{\R\bar\R X \ar[r]^-{\R\bar\phi'_X} & \R \bar\L X^{m'} \ar[r]^-{\phi'_{\bar\L X^{m'}}} & \L\bar\L X^{m+m'}} ,
&&
\xymatrix{
\L\bar\L X^{m+m'} \ar[r]^-{\psi'_{\bar\L X^{m'}}} & \R\bar\L X^{m'} \ar[r]^-{\R\bar\psi'_X} & \R\bar\R X};
\end{eqnarray*}
which are easily verified to satisfy the required properties.
\end{proof}

The notion of a \emph{left quasi-Frobenius pair of functors}
was introduced by Guo in \cite{guo}, where he proved that a ring extension
$\varphi: R\rightarrow S$ is a \emph{left quasi-Frobenius extension}
if and only if $(\varphi_*, -\otimes_RS)$  is a \emph{left
quasi-Frobenius pair of functors}.
\\
In general, for categories $\mathcal{A}$ and $\mathcal{B}$ with
finite direct sums,  the pair of functors
$(\functor{F},\functor{L})$ is called a \emph{left quasi-Frobenius
pair of functors} if $\functor{F}:\mathcal{A}\rightarrow
\mathcal{B}$ is a right adjoint of
$\functor{L}:\mathcal{B}\rightarrow \mathcal{A}$ and for some
positive integer $n$, there are natural transformations
$$\alpha: 1_{\mathcal{A}}\rightarrow (\functor{LF})^n  \, \,
\mbox{   and } \, \,   \overline{\alpha}: (\functor{FL})^n
\rightarrow
 1_{\mathcal{B}}$$
such that $$\overline{\alpha}_{\functor{F}(X)}\circ
\functor{F}(\alpha_X) = 1_{\functor{F}(X)}
$$
for all $X\in \mathcal{A}$.

Recall that if the functor
$\functor{R}:\mathcal{B}\rightarrow \mathcal{A}$ is a right adjoint
to $\functor{F}$, then the unit $\overline{\eta}:
1_{\mathcal{A}}\rightarrow \functor{RF}$ and the counit
$\overline{\rho}: \functor{FR}\rightarrow 1_{\mathcal{B}}$ of the adjunction satisfy
the identities $\overline{\rho}_{\functor{F}(X)}\circ
\functor{F}(\overline{\eta}_X) = 1_{\functor{F}(X)}$ and
$\functor{R}(\overline{\rho}_Y)\circ
\overline{\eta}_{\functor{R}(Y)} = 1_{\functor{R}(Y)},$ for all
$X\in \mathcal{A}$ and $Y\in \mathcal{B}$.

The next proposition implies, in particular, that if $\Gamma
=(\functor{L,F,R})$ is a quasi-Frobenius triple of functors, then
$(\functor{F,L})$ is a left quasi-Frobenius pairs. For a more complete treatment about the connection between quasi-Frobenius triples and (left) quasi-Frobenius pairs, we refer to \cite{Ver:QcF}.
\begin{proposition}  If $\Gamma =(\functor{L,F,R})$ is a quasi-Frobenius
triple of functors for $\mathcal{A}$ and $\mathcal{B}$, then
$(\functor{F,L})$ and $(\functor{R,F})$ are left quasi-Frobenius
pairs.
\end{proposition}
\begin{proof} Assume $\functor{R}|\functor{L}$. Then there
exist  morphisms
\[\xymatrix{\functor{R}\ar[r]^-{\phi'}&  \functor{L}^n
 \ar[r]^-{\psi'}& \functor{R} }\]
such that $\psi'\circ\phi' = 1_{\R}.$

We define  the functorial morphism $\alpha:
1_{\mathcal{A}}\rightarrow (\functor{LF})^n$ by the composition of
morphisms
\[\xymatrix{X  \ar[r]^-{\overline{\eta}_X}&\functor{RF}(X)\ar[r]^-{\phi'_{\functor{F}(X)}}&  \functor{L}^n
 (\functor{F}(X))\cong (\functor{LF})^n
 (X) }\]
for every $X\in \mathcal{A}$. Similarly, for any $Y\in \mathcal{B},$
the composition
\[\xymatrix{(\functor{FL})^n
 (Y)\cong  \functor{F}(\functor{L}^n
 (Y)) \ar[r]^-{\functor{F}(\psi'_Y)}&\functor{FR}(Y)\ar[r]^-{\overline{\rho}_Y}&  Y }\]
defines the functorial morphism $\overline{\alpha}:
(\functor{FL})^n \rightarrow 1_{\mathcal{B}} $. Then
$$
\begin{array}{ccl} \overline{\alpha}_{\functor{F}(X)}\circ \functor{F}(\alpha_X) & = &
\overline{\rho}_{\functor{F}(X)}\circ
\functor{F}(\psi'_{\functor{F}(X)}\circ \phi'_{\functor{F}(X)}\circ
\overline{\eta}_X)
          \\[+1mm]
& = & \overline{\rho}_{\functor{F}(X)}\circ
\functor{F}(1_{\functor{RF}(X)}\circ \overline{\eta}_X) \\[+1mm]
& = & \overline{\rho}_{\functor{F}(X)}\circ\functor{F}(
\overline{\eta}_X) = 1_{\functor{F}(X)}.
\end{array}$$
This means that $(\functor{F,L})$ is a left quasi-Frobenius pair.
Likewise, using that $\functor{L}|\functor{R}$, we obtain the other
afirmation.
\end{proof}

\section{Quasi-Frobenius functors between module categories}\selabel{modules}

\begin{lemma}\lelabel{objdivfundiv}
Let $M$ and $N$ be $(R,S)$-bimodules, and suppose that $M~|~N$. If $F:{_R\Mm_S}\to{_{R'}\Mm_{S'}}$ preserves finite direct sums, then $F(M)~|~F(N)$ as $(R',S')$-bimodules. Consequently, if $M\sim N$ as $(R,S)$-bimodules, then $F(M)\sim F(N)$ as $(R',S')$-bimodules.
\end{lemma}

\begin{proof}
Suppose that $M~|~N$, by \leref{splitting} this is true if and only if there exists a $P\in{_R\Mm_S}$ such that the following sequence splits in ${_R\Mm_S}$.
$$\xymatrix{0\ar[r] & M \ar[r] & N^{(n)} \ar[r] & P \ar[r] & 0}.$$
Recall that any functor preserves split exact sequences and $F$ preserves finite direct sums. Therefore, if we apply the functor $F$ to the above sequence, then we obtain the split exact sequence
$$\xymatrix{0\ar[r] & F(M) \ar[r] & F(N)^{(n)} \ar[r] & F(P) \ar[r] & 0}$$
in ${_{R'}\Mm_{S'}}$. Hence $F(M)~|~F(N)$.
\end{proof}

\begin{lemma}\label{simi}
Let $M$ and $N$ be $(R,S)$-bimodules. Then
\begin{enumerate}[(i)]
\item
$_RM_S~|~{_RN_S}$ if and only if $M\ot_S-~|~N\ot_S-$;
\item
$_RM_S\sim \, _RN_S$ if and only if
the tensor functors $M\otimes_{S}- $ and $N\otimes_{S}-$ are
similar.
\end{enumerate}
\end{lemma}
\begin{proof} We only prove part (i). Suppose first that $M\ot_S- ~|~ N\ot_S-$, then the statement follows by \leref{fundivobjdiv}. Conversely, if $M~|~N$, then also $M\ot_SP~|~N\ot_SP$ for all $P$, because the tensor product functor $-\ot_SP$ preserves direct sums and therefore we can apply \leref{objdivfundiv}. Naturality of this divisibility is easily checked.
\end{proof}

\begin{lemma}\lelabel{divfgp}
Let $M$ and $N$ be $(R,S)$-bimodules. Suppose that $M$ divides $N$. Then
\begin{enumerate}[(i)]
\item if $N$ is finitely generated and projective as a right $S$-module, then $M$ is finitely generated and projective as a right  $S$-module;
\item if $N$ is finitely generated and projective as a left $R$-module, then $M$ is finitely generated and projective as a left $R$-module.
\end{enumerate}
\end{lemma}

\begin{proof}
We only prove part (i), the proof of the second part is similar. Take a finite dual basis $\{e_k,f_k\}\in N\times \Hom_S(N,S_S)$ for $M$ as a right $S$-module. Since $M$ divides $N$, there exist morphisms
$$\phi:M\to N^{(n)}, \qquad \psi:N^{(n)}\to M,$$
such that $\psi\circ\phi= 1_M$. Denote by $\pi_i:N^{(n)}\to N$ and $\iota_i:N\to N^{(n)}$ the canonical projection and injection on the $i$-th component of the direct sum.
Take any $m\in M$, then
\begin{eqnarray*}
m&=&\psi\circ \phi(m)= \psi(\sum_i\iota_i\pi_i(\phi(m)))=\sum_i\psi\iota_i( \pi_i\phi(m))\\
&=&\sum_i\psi\iota_i(e_k f_k(\pi_i\phi(m)))=\sum_i\psi\iota_i(e_k) f_k(\pi_i\phi(m))
\end{eqnarray*}
Hence we obtain a dual basis $\{\psi\circ\iota_i(e_k), f_k\circ\pi_i\circ\phi \}\in M\times\Hom_S(M,{S_S})$.
\end{proof}

Quasi-Frobenius functors between module categories have already been considered by K.\ Morita in \cite{Morita}. The following theorem extends the characterization given in that paper.
\begin{theorem}\label{modules}\thlabel{Morita}
 For functors $\functor{F}:
\lmod{R}\rightarrow \lmod{S}$ and $\functor{L,R}:
\lmod{S}\rightarrow \lmod{R},$ the following statements are
equivalent.
\begin{enumerate}[(a)]
\item $(\functor{L,F,R})$ is a quasi-Frobenius
triple;
\item There exist bimodules $_SM_R, \, _RN_S$ and $_R\overline{N}_S
$ with the following properties:
\begin{enumerate}[(i)]
\item $_SM_R$ is finitely generated and
projective on both sides;
\item ${_RN_S}\cong \, (M_R)^* $ and $_R\overline{N}_S \cong \, (_SM)^*$ as bimodules;
\item $\functor{F}\cong M\otimes_R -, \, \, \functor{L}\cong N\otimes_S - $ and \, $\functor{R}\cong \overline{N} \otimes_S -$;
\item $_RN_S \sim \,_R\overline{N}_S$.
\end{enumerate}
\item There exist bimodules satisfying (b)(iii) and (iv), as well as
\begin{enumerate}[(i)']
\item $_RN_S$ is finitely generated and projective on both sides
\item $M\cong (_RN)^*$ and $\ol N\cong (_S(_RN)^*)^*$;
\end{enumerate}
\item There exist bimodules satisfying (b)(iii) and (iv), as well as
\begin{enumerate}[(i)'']
\item $_R\ol N_S$ is finitely generated and projective on both sides
\item $M\cong (\ol N_S)^*$ and $N\cong ((\ol N_S)^*\; _R)^*$;
\end{enumerate}
\end{enumerate}
\end{theorem}
\begin{proof}
The equivalence $(a)\Leftrightarrow(b)$ is proven in {\cite[Theorem 4.1]{Morita}}. The equivalence with $(c)$ and $(d)$ follows immediately from \leref{divfgp} and the characterization of finitely generated and projective modules in terms of adjoint functors that can be found for example in \cite[Theorem 3.1]{Morita}.
\end{proof}

The following proposition shows that there is a duality between the category of quasi-Frobenius functors for ${_R\Mm}$ and ${_S\Mm}$ on one hand and the category of quasi-Frobenius functors for ${\Mm_S}$ and $\Mm_R$ on the other hand.

\begin{proposition}\prlabel{symmetry}
For each quasi-Frobenius triple $(\functor{L,F,R})$ for $\lmod{R}$
and $\lmod{S}$ there is a quasi-Frobenius triple
$(\functor{\overline{L},\overline{F},\overline{R}})$ for $\rmod{S}$
and $\rmod{R}$ such that the correspondence
$(\functor{L,F,R})\mapsto
(\functor{\overline{L},\overline{F},\overline{R}})$ between
quasi-Frobenius triples is bijective up to natural isomorphism.
\end{proposition}

\begin{proof}
Let ${_SM_R}$, ${_RN_S}$ and ${_R\ol{N}_S}$ be bimodules associated to the triple $(\L,\F,\R)$ as in part (b) of \thref{Morita}. Then we define $\ol{\L}=-\ot_R\ol{N}$, $\ol{\F}=-\ot_SM$ and $\ol{\R}=-\ot_R N$. By the left-right dual of \thref{Morita}, we find that $(\ol{\L},\ol{\F},\ol{\R})$ is indeed a quasi-Frobenius triple and this correspondence is clearly bijective.
\end{proof}

Recall from \cite{A-F} that an $(R,S)$-bimodule $M$ is called \emph{Frobenius} if $_RM$ and
$M_S$ are finitely generated and projective and $(_RM)^*\cong \,
(M_S)^*$ as $(S,R)$-bimodules.
Motivated by the previous theorem, we propose in this note the following generalization of Frobenius bimodules. It should be remarked however, that this notion differs from the more classical definition of a quasi-Frobenius bimodule (see \cite{Azu:QF}).
\begin{definition}\label{qF-bimodule}\delabel{qfbim} An  $(R,S)$-bimodule $M$ is
said to be \emph{quasi-Frobenius bimodule}, if both $_RM$ and $M_S$
are finitely generated projective and $(_RM)^*\sim \, (M_S)^*$ as
$(S,R)$-bimodules.
\end{definition}

\thref{Morita} and \prref{symmetry} immediately lead to the following.
\begin{proposition}\label{bimo-funtor} For any $(R,S)$-bimodule $M$ the following
assertions are equivalent.
 \begin{enumerate}[(i)]
 \item $_RM_S$ is a quasi-Frobenius bimodule;
 \item $M\otimes_S-:{_S\Mm}\to {_R\Mm}$ is a quasi-Frobenius functor;
 \item[(ii)'] there is a quasi-Frobenius functor $\F:{_S\Mm}\to {_R\Mm}$ such that $\F(S)\cong M$ as $(R,S)$-bimodule;
 \item $-\ot_RM:\Mm_R\to\Mm_S$ is a quasi-Frobenius functor;
\item[(iii)'] there is a quasi-Frobenius functor $\bar{\F}:{\Mm_R}\to {\Mm_S}$ such that $\bar{\F}(R)\cong M$ as $(R,S)$-bimodule.
 \end{enumerate}
\end{proposition}

The following theorem shows that a quasi-Frobenius triple between module categories can `shifted'.

\begin{theorem}\thlabel{sequence}
Consider an adjoint triple $(\L,\F,\R)$ for ${_R\Mm}$ and ${_S\Mm}$. Then the following statements are equivalent
\begin{enumerate}[(i)]
\item $(\L,\F,\R)$ is a quasi-Frobenius triple, i.e.\ $\F$ is a quasi-Frobenius functor;
\item there exists a functor $\F_1 :{_R\Mm}\to {_S\Mm}$ such that $(\F,\R,\F_1)$ is a quasi-Frobenius triple for ${_S\Mm}$ and ${_R\Mm}$, i.e. $\R$ is a quasi-Frobenius functor;
\item there exists a functor $\F_{-1}:{_R\Mm}\to {_S\Mm}$ such that $(\F_{-1},\L,\F)$ is a quasi-Frobenius triple for ${_S\Mm}$ and ${_R\Mm}$, i.e.\ $\L$ is a quasi-Frobenius functor.
\end{enumerate}

\end{theorem}
\begin{proof}
We only prove $(i)\Rightarrow (ii)$ and $(i)\Rightarrow(iii)$. The converse implications follow from this by symmetry arguments.

Let ${_SM_R}$, ${_RN_S}$ and ${_R\ol N_S}$ be bimodules associated to the triple $(\L,\F,\R)$ as in \thref{Morita} (d). Now put $\F_1=(_R{\ol N})^*\ot_R-$. Since $\ol N$ is finitely generated and projective as a left $R$-module, $\F_1$ is a right adjoint for $\R$. Moreover, we know that $N\sim \ol N$, therefore $M\cong {\Hom_R}(N,{_RR}) \sim \Hom_R(\ol N,{_RR})$. So by applying again \thref{Morita} ($(b)\Rightarrow(a)$), we find that $(\F,\R,\F_1)$ is a quasi-Frobenius triple.

Similarly, we obtain the quasi-Frobenius triple $(\F_{-1},\L,\F)$, by putting $\F_{-1}=(N_S)^*\ot_R-$.
\end{proof}

\begin{remark}\relabel{sequence}
Let $M\in{_S\Mm_R}$ be a quasi-Frobenius bimodule. Then we obtain a sequence of functors
\[
\hspace{-.8cm}
\xymatrix{
\ldots \to {_R\Mm} \ar[rr]^-{((M_R)^*\,_S)^*\ot_R-} &&
{_S\Mm} \ar[rr]^-{ (M_R)^*\ot_S-} &&
{_R\Mm} \ar[r]^-{ M\ot_R-}  &
{_S\Mm} \ar[rr]^-{ (_SM)^*\ot_S-} &&
{_R\Mm} \ar[rr]^-{ (_R(_SM)^*)^*\ot_R-} &&
{_S\Mm} \to \ldots
}
\]
such that each three subsequent functors form a quasi-Frobenius triple of functors. In particular, a bimodule is quasi-Frobenius if and only if any of it's duals appearing in this sequence is quasi-Frobenius.
\end{remark}

If $_RM_S$ and $_SN_T$ are bimodules, then $M\otimes_SN$ receives the natural $(R,T)$-bimodule structure
by putting $r(m\ot_S n)t=rm\ot_S nt$ for all $r\in R$, $t\in T$ and $m\ot_S n\in M\ot_S N$.
\begin{proposition}\label{product} Suppose $_RM_S$ and $_SN_T$ are quasi-Frobenius bimodules. Then $_R(M\otimes_SN)_T$ is  a
quasi-Frobenius $(R,T)$-bimodule.
\end{proposition}
\begin{proof}
This follows directly from Proposition \ref{properties} and Lemma \ref{simi}.
\end{proof}

\section{Quasi-Frobenius ring extensions}

Let $\varphi:R\to S$ be a ring morphism. Then $\varphi$ induces an $R$-bimodule structure on $S$. Recall from \cite{muller} that $\varphi$ is called a quasi-Frobenius ring extension if and only if $S$ is finitely generated and projective as a left $R$-module and the $(S,R)$-bimodules $S$ and $\Hom_R({_RS},{_RR})$ are similar. From this definition, it is clear that $\varphi$ is a quasi-Frobenius ring extension if and only if ${_RS_S}$ is a quasi-Frobenius bimodule.

\begin{remark}
It is obvious that the class of quasi-Frobenius bimodule contains Frobenius bimodules.
An example of finite-dimensional quasi-Frobenius algebra which is not a Frobenius is given in
\cite{nakayama}. Consequently, a quasi-Frobenius bimodule need not be a Frobenius in general.
\end{remark}

We can associate to the ring morphism $\varphi:R \to S$ the adjoint triple of
functors $\Gamma =(S\otimes_R-, {^l\varphi_*}, \hom{R}{_RS_S}{-})$
where ${^l\varphi_*}: \lmod{S}\rightarrow \lmod{R}$ is the restriction
of scalars functor. By symmetry, we can also consider  the
restriction functor $\varphi_*^r: \rmod{S}\rightarrow \rmod{R}$ and
its adjoint functors $-\otimes_RS$ and $\hom{R}{S_R}{-}.$
By applying the results of the previous section, we obtain the following characterization of quasi-Frobenius ring extensions.

\begin{corollary}\colabel{frobext} Let $\varphi:R\rightarrow S$ be a ring extension.
Then the following assertions are  equivalent
\begin{enumerate}[(i)]
\item $\varphi$ is a quasi-Frobenius extension;
\item $_RS_S$ is a quasi-Frobenius bimodule;
\item[(ii)'] $_SS_R$ is a quasi-Frobenius bimodule;
\item $\Gamma = (S\otimes_R-, \, ^l\varphi_*, \, \hom{R}{_RS}{-})$ is a
quasi-Frobenius triple for $\lmod{R}$ and $\lmod{S}$, i.e.\ ${^l\varphi_*}: \lmod{S}\rightarrow \lmod{R}$ is a quasi-Frobenius functor;
\item[(iii)'] $\overline{\Gamma} = (-\otimes_RS, \, \varphi_*^r, \,
\hom{R}{S_R}{-})$ is also a quasi-Frobenius triple for $\rmod{R}$
and $\rmod{S}$, i.e.\ $\varphi_*^r:\Mm_S\to\Mm_R$ is a quasi-Frobenius functor;
\item $-\ot_RS:\Mm_R\to \Mm_S$ is a quasi-Frobenius functor;
\item[(iv)'] $S\ot_R-:{_R\Mm}\to {_S\Mm}$ is a quasi-Frobenius functor.
\end{enumerate}
\end{corollary}

\begin{proof} The equivalences follow now easily by combining the results of the previous section. Just observe that $-\otimes_SS \simeq \varphi_*^r:\Mm_S\to\Mm_R$ and $S\ot_S-\simeq {^l\varphi_*}:{_S\Mm}\to {_R\Mm}$.
\end{proof}

\begin{corollary}
Let $\alpha:R\to S$ and $\beta:S\to T$ be two ring morphisms. Suppose that $\beta$ is a quasi-Frobenius extension. Then $\alpha$ is a quasi-Frobenius extension if and only if $\beta\circ\alpha$ is a quasi-Frobenius extension.
\end{corollary}

\begin{proof}
Consider the following diagram of functors
\[
\xymatrix{
\Mm_R \ar@<-.5ex>[ddrr]_-{-\ot_RT} \ar@<.5ex>[rrrr]^-{-\ot_RS} &&&& \Mm^\dd \ar@<.5ex>[ddll]^-{-\ot_ST} \ar@<.5ex>[llll]^-{U_\alpha}
\\  \\
&& \Mm_T \ar@<-.5ex>[uull]_-{U_{\beta\circ\alpha}} \ar@<.5ex>[uurr]^-{U_\beta}
}
\]
where the indexed functors $U_{-}$ are the obvious forgetful functors and the inner and outer triangles are commutative. We know that the extensions are quasi-Frobenius if either the forgetful or the induction functors are quasi-Frobenius functors (see \coref{frobext}).
The statements follow now directly form Lemma \ref{properties} (b). Remark that the `only if' part is in fact a special instance of Proposition \ref{product}.
\end{proof}

\section{Quasi-Frobenius functors in graded rings}\selabel{graded}

Let $G$ be a group with neutral element $e$. A ring $R$ is said to
be $G$-graded if there is a family $\{R_x ; x\in G\}$ of additive
subgroups of $R$ such that $R= \bigoplus_{x\in G}R_x$, and the
multiplication in $R$ is such that, for all $x$ and $y$ in $G$,
$R_xR_y\subseteq R_{xy}$. Similarly, a left $R$-module $M$ is graded
by $G$ if there is a family $\{M_x ; x\in G\}$ of additive subgroups
of $M$ such that $M= \bigoplus_{x\in G}M_x$, and
 for all $x$ and $y$ in $G$,
$R_xM_y\subseteq M_{xy}$.
  We will denote by $R$-gr the category of all $G$-graded left
$R$-modules over the unital group-graded ring $R$.

It is well known (see e.g. \cite{NO:2004}) that associated to the ring homomorphism
$\varphi : R_e\rightarrow R$ can associate two functors
$$Ind(-): {_{R_e}\Mm} \to R\mbox{-gr} \quad {\rm and} \quad
(-)_e: R\mbox{-gr} \to {_{R_e}\Mm}.$$
Here the functor $(-)_e$, called the \emph{restriction at $e$}, is given by
$M\mapsto M_e$, for every left graded $R$-module
 $M=\bigoplus _{x\in G}M_x$ and $Ind(-)$. The functor $Ind(-)$, called the \emph{induction} functor
is given by $Ind(N)=R\otimes_{R_e}N$, for every left $R_e$-module $N$, where the grading on $Ind(N)$ is defined by putting
 $$(Ind(N))_y = R_y\otimes _{R_e}N $$
 for every $y\in G$.
It was shown in
\cite{NO:2004} that the functor $Ind(-)$ is a left adjoint of the
functor  $(-)_e$ and the unit of the adjunction
 $\eta: 1_{\lmod{R_e}}\rightarrow
(-)_e\circ Ind(-)$ is a functorial isomorphism.

The functor $(-)_e$ has also a right adjoint called the \emph{e-th
coinduced} functor $$Coind(-): \lmod{R_e}\rightarrow R\mbox{-gr},$$
where for every left $R_e$-module $N$, $Coind(N)=\bigoplus_{y
in G} Coind(N)_y$ is the left graded $R$-module defined by
$$Coind(N)_y = \{ f\in \hom{R_e}{R}{N} \, | \, \, f(R_x) =0, \, \forall x\neq y^{-1} \}\cong \Hom_{R_e}(R_{y^{-1}},N)$$
Moreover, the counity of this adjunction  $\tau: (-)_e\circ
Coind(-)\rightarrow 1_{\lmod{R_e}}$ is a functorial isomorphism.
From our definitions, it follows that $(-)_e$ is a quasi-Frobenius functor if and only if
$(Ind(-), (-)_e,Coind(-))$ is a quasi-Frobenius triple of functors.

Recall that for any $M\in R\mbox{-gr}$ and $x\in G$, we define the $x$-suspension $M(x)$ of $M$ as the graded
$R$-module obtained from $M$ by putting $M(x)_y=M_{yx}$ for all $y\in G$.

\begin{theorem} Let $R$ be a $G$-graded ring. The following assertions are equivalent.
\begin{enumerate}[(i)]
\item $(-)_e$ is a quasi-Frobenius functor;
\item $Ind(-) \sim Coind(-)$;
 \item
$\forall \, x \in G, \, R_x$ is finitely generated and projective in
$\lmod{R_e}$ and $R\sim Coind(R_e)$.
\end{enumerate}
\end{theorem}
\begin{proof}
$\ul{(i)\Leftrightarrow (ii)}$ is clear.
\\
$\ul{(ii)\Rightarrow (iii)}.$  Assume that $Ind(-)\sim Coind(-)$. Then
$Ind(R_e)\sim Coind(R_e).$ But $Ind(R_e)\cong R$ which implies that
$R\sim Coind(R_e)$ as $(R,R_e)$-bimodules.
Now consider the $x$-suspended objects $R(x)\in R\mbox{-gr}$ for all $x\in G$. These are finitely generated and projective in $R\mbox{-gr}$, which can easily be seen from the fact that the forgetful functor $R\mbox{-gr}\to {_R\Mm}$ reflects finitely generated and projective objects.
Furthermore, by Lemma \ref{properties} (a), the functor $(-)_e$ preserves finitely generated and projective objects. Hence $(R(x))_e=R_x$ is a finitely generated and projective left $R_e$-module for all $x\in G$.
\\
$\ul{(iii)\Rightarrow (ii)}.$ Assume that $_RR_{R_e} \, | \, Coind(R_e)$.
Then there exist morphisms in $\lmod{R_e}$
\begin{equation}\label{dos} \xymatrix{R\ar[r]^-{f}&
Coind(R_e)^n
 \ar[r]^-{g}& R }
\end{equation}
with $g\circ f =1_R$. For any left $R_e$-module $X$, apply the
functor $-\otimes _{R_e} X$ to (\ref{dos}) and we obtain
\begin{equation}\label{tres}
 \xymatrix{R\otimes_{R_e}X \ar[r]^-{f\otimes X}&
Coind(R_e)^n\otimes_{R_e}X\cong (Coind(R_e)\otimes_{R_e}X)^n
 \ar[r]^-{g\otimes X}& R \otimes_{R_e}X}
 \end{equation}
By assumption $R_x$ is finitely generated and projective as $R_e$-module, whence
$$\hom{R_e}{R_x}{R_e}\otimes_{R_e}X\cong \hom{R_e}{R_x}{X}.$$ In
particular, $Coind(R_e)\otimes_{R_e}X \cong Coind(X).$ Then the
sequence (\ref{tres}) is given by
\[
\xymatrix{R\otimes_{R_e}X \ar[r]^-{f\otimes X}& Coind(X)^n
 \ar[r]^-{g\otimes X}& R \otimes_{R_e}X}
 \]
Since  $(g\otimes X)\circ (f\otimes X)= (g\circ f)\otimes X =
1_R\otimes X,$ this implies that $Ind(-) \, | \, Coind(-).$
 Analogously, we can prove that $Coind(-) \, | \, Ind(-)$. Therefore, $Ind(-)\sim
Coind(-)$.
\end{proof}
\begin{remark}\relabel{qf=f}
Let $R= \bigoplus_{x\in G}R_x$ be a $k$-algebra graded by a group
$G$. We consider the forgetful functor $U:R\mbox{-gr}\rightarrow
\lmod{R}$, where $R$-gr is the category of $G$-graded modules. It is
well know that $U$ has a right adjoint functor $F:\lmod{R}
\rightarrow R$-gr. If $U$ is a quasi-Frobenius functor, then  $U$
commutes with direct products and by \cite[Corollary 4.4]{CGN}, $G$
is finite. This implies that $U$ is a Frobenius functor (see
\cite[Proposition 2.5]{RDNO}).
\end{remark}

\section{Quasi-Frobenius functors between  categories of comodules over corings}\selabel{comodules}

Troughout this section, let $\cc$ be an $A$-coring and $\dd$ a $B$-coring, where $A$ and $B$ are $k$-algebras over the commutative ring $k$.

\begin{lemma}\lelabel{functorcomod}
Suppose that $_A\mathfrak{C}$ and $_B\mathfrak{D}$ are flat and $N,
\overline{N}\in \, ^{\mathfrak{D}}\mathcal{M}^{\mathfrak{C}}$. Then
$_{\mathfrak{D}}N_{\mathfrak{C}}\sim \,
_{\mathfrak{D}}\overline{N}_{\mathfrak{C}}$ if and only if $\,
-\cotensor{\mathfrak{D}}N \sim \,
-\cotensor{\mathfrak{D}}\overline{N}.$
\end{lemma}
\begin{proof} If $-\cotensor{\mathfrak{D}}N \sim \, -\cotensor{\mathfrak{D}}\overline{N}$, then
$\mathfrak{D}\cotensor{\mathfrak{D}}N \sim \,
\mathfrak{D}\cotensor{\mathfrak{D}}\overline{N}$ and hence
$_{\mathfrak{D}}N_{\mathfrak{C}}\sim \,
_{\mathfrak{D}}\overline{N}_{\mathfrak{C}}$. Assume now that
$_{\mathfrak{D}}N_{\mathfrak{C}}\mid \,
_{\mathfrak{D}}\overline{N}_{\mathfrak{C}}$. This  condition
establishes  that for some positive integer $n$ there are bicomodule
morphisms
\[\xymatrix{N\ar[r]^-{f}&  \overline{N}^{\; n}
 \ar[r]^-{g}& N }\]
such that $g\circ f = 1_{N}.$ For any right $\mathfrak{D}$-comodule
$X$ we apply the cotensor functor $X\cotensor{\mathfrak{D}}-$ to the
above sequence and we obtain
\[\xymatrix{X\cotensor{\mathfrak{D}}N\ar[r]^-{1_X\cotensor{}f}&  (X\cotensor{\mathfrak{D}}\overline{N})^n
 \ar[r]^-{1_X\cotensor{}g}& X\cotensor{\mathfrak{D}}N }\]
Clearly, $(1_X\cotensor{}g)\circ (1_{X}\cotensor{}f) =
1_X\cotensor{}(g\circ f)= 1_{X\cotensor{\mathfrak{D}}N}.$ This
implies that $-\cotensor{\mathfrak{D}}N \mid
-\cotensor{\mathfrak{D}}\overline{N}.$ Analogously, from $
_{\mathfrak{D}}\overline{N}_{\mathfrak{C}} \mid \,
_{\mathfrak{D}}N_{\mathfrak{C}}$ we get
$-\cotensor{\mathfrak{D}}\overline{N} \mid
-\cotensor{\mathfrak{D}}N.$  Thus $-\cotensor{\mathfrak{D}}N \sim
-\cotensor{\mathfrak{D}}\overline{N}.$
\end{proof}

We recall from \cite{Ka} that a bicomodule $X \in {^\cc\Mm^\dd}$ is called an $(A,\dd)$-\emph{injector} if the functor
$-\otimes_AX:\Mm_A\to \Mm^\dd$
preserves injective objects. Furthermore, $X$ is said to be \emph{$(A,\dd)$-quasi-finite} if the functor $-\ot_AX:\Mm_A \to \Mm^\dd$ has a left adjoint, which will in this case be denoted by ${\rm h}_\dd(X,-):\Mm^\dd\to \Mm_A$ and is called the \emph{Cohom-functor}. In case $\dd$ is flat as a left $B$-module than the Cohom-functor can be restricted to a functor $\Mm^\dd\to \Mm^\cc$, being a left adjoint for $-\cot_\dd X:\Mm^\dd\to\Mm^\cc$. In fact the existence of this restricted adjoint is in this case equivalent to $X$ being quasi-finite, see \cite[Proposition 4.2]{gomez} or \cite[23.6]{BrzWis:book}.

\begin{theorem}\label{main}
Suppose that  $_A\mathfrak{C}$ and $_B\mathfrak{D}$ are flat. For
$k$-linear functors $\functor{F}: \rcomod{\mathfrak{C}}\rightarrow
\rcomod{\mathfrak{D}}$ and $\functor{L,R}:
\rcomod{\mathfrak{D}}\rightarrow \rcomod{\mathfrak{C}},$ the
following statements are equivalent.
\begin{enumerate}[(a)]
\item $\Gamma = (\functor{L,F,R})$ is a quasi-Frobenius triple.
\item There exist bicomodules ${_\cc M_\dd}, \, {_\dd N_\cc}$ and
${_\dd\ol{N}_\cc} $ with the following
properties.
\begin{enumerate}[(i)]
\item
$\functor{L}\cong-\cot_\dd N$,
$\functor{F}\cong -\cot_\cc M$,
and $\functor{R}\cong -\cotensor{\mathfrak{D}}\overline{N}$.
\item ${_AM_\dd}$ and ${_B\overline{N}_\cc}$ are quasi-finite injectors.
\item $\cohom{\mathfrak{C}}{\overline{N}}{\mathfrak{C}}\cong \, _{\mathfrak{C}}M_{\mathfrak{D}},$
$\cohom{\mathfrak{D}}{M}{\mathfrak{D}}\cong \,
_{\mathfrak{D}}N_{\mathfrak{C}}$;
\item
$N\sim \ol{N}$ as $(\dd,\cc)$-bicomodules.
\end{enumerate} \end{enumerate}
\end{theorem}
\begin{proof}$\ul{(a) \Rightarrow (b)}.$
From Lemma \ref{properties}, we know that the functors $\L$, $\F$ and $\R$ are exact and preserve all limits and colimits. Therefore, it follows that there exist bimodules $M$, $N$ and  $\ol{N}$ as in part $(i)$ of statement $(b)$ (see e.g. \cite[Theorem 3.5]{gomez} or \cite[23.1]{BrzWis:book}.
Assertions $(ii)$ and $(iii)$ follow now as an immediate consequence of
\cite[Proposition 2.9]{Za}.
 Finally, by \leref{functorcomod} $\functor{L}\sim
 \functor{R}$ implies that $N\sim \ol{N}$ as $(\dd,\cc)$-bicomodules.
 \\
 $\ul{(b) \Rightarrow (a)}.$ Assume that there exist bicomodules $_{\mathfrak{C}}M_{\mathfrak{D}}, \, _{\mathfrak{D}}N_{\mathfrak{C}}$ and
$_{\mathfrak{D}}\overline{N}_{\mathfrak{C}} $ satisfying all conditions of part $(b)$.
Combining conditions $(i)$ and $(iv)$ with \leref{functorcomod}, we obtain that $\L\sim\R$.
Using the property of quasi-finite comodules recalled before this theorem,
we know that $\F\cong -\cot_\cc M$ has a left adjoint $\cohom{\cc}{M}{-}$ and $\R\cong-\cot_\dd\ol{N}$ has a left adjoint $\cohom{\dd}{\ol{N}}{-}$. Using a property of quasi finite injectors (see \cite[Corollary 3.12]{Ka}) in the last isomorphism of the next computation, we find that for all $X\in\Mm^\dd$
$$\L(X)\cong X\cot_\dd\ol{N}\cong X\cot_\dd\cohom{\dd}{M}{\dd}\cong \cohom{\dd}{M}{X}.$$
Hence, $\L$ is a left adjoint for $\F$. Similarly, we find that $\F$ is a left adjoint for $\R$, and therefore
$\Gamma$ is a quasi-Frobenius
triple.
\end{proof}

\section{Quasi-Frobenius coring homomorphisms} \selabel{corings}

Following \cite{gomez}, a coring homomorphism from the $A$-coring
$\cc$ to the $B$-coring $\dd$ is a pair
$(\varphi, \rho)$, where $\rho:A\rightarrow B$ is a homomorphism of
$k$-algebras and $\varphi:\cc\to \dd$ is a
homomorphism of $A$-bimodules such that
$$\varepsilon_\dd \circ \varphi =\rho \circ \varepsilon_\cc \mbox{  and   }
\Delta_\dd\circ \varphi = \omega_{\dd,\dd}\circ (\varphi\otimes_A\varphi)\circ \Delta_{\cc},$$
where $\omega_{\dd,\dd}:\dd\ot_A\dd\to\dd\ot_B\dd$ is the canonical map induced by
$\rho:A\to B$. The functor $-\ot_AB:\Mm^\cc\to\Mm^\dd$ has a right adjoint
$-\cot_\dd(B\ot_A\cc):\Mm^\dd\to\Mm^\cc$, hence $B\ot_A\cc$ is $(B,\cc)$ quasi-finite  \cite[Proposition 5.4]{gomez}. From now on, suppose that ${_A\cc}$ and ${_B\dd}$ are flat.
If $(\cc\ot_AB)$ is $(A,\dd)$ quasi-finite, then the functor $\cohom{\dd}{\cc\ot_AB}{-}:\Mm^\dd\to\Mm^\cc$ is a left adjoint to
$-\cot_\cc(\cc\ot_AB)\cong (-\cot_\cc\cc)\ot_AB\cong -\otimes_AB$.
In this case we have an adjoint triple of functors $$\Gamma = (\cohom{\dd}{\cc\ot_AB}{-}, -\ot_AB,-\cot_\cc(B\ot_A\cc))$$ between the Grothendieck categories $\Mm^\cc$ and $\Mm^\dd$.
Moreover, $(B\ot_A\cc)$ is a $(B,\cc)$ (quasi-finite) injector because the functor $-\cot_\dd(B\ot_A\cc)$ is right adjoint to the exact functor $-\ot_AB$. These observations, in combination with Theorem \ref{main}, lead to the following.

\begin{theorem}\label{qF-morphism}
Let $\cc$ be an $A$-coring and $\dd$ a $B$-coring such that ${_A\cc}$ and ${_B\dd}$ are flat.
Consider a homomorphism of corings $(\varphi, \rho):\cc\to\dd$ such that $B\ot_A\cc$ is $(B,\cc)$ quasi-finite. Then the following statements are equivalent.
\begin{enumerate}[(i)]
\item $-\otimes_AB:\mathcal{M}^{\mathfrak{C}}\rightarrow
\mathcal{M}^{\mathfrak{D}}$ is a quasi-Frobenius functor;
\item $\mathfrak{C}\otimes_AB $ is an $(A,\dd)$ quasi-finite injector and $\cohom{\dd}{\cc\ot_AB}{\dd}\sim B\ot_A\cc$ as $(\dd,\cc)$-bicomodules.
\end{enumerate}
\end{theorem}

\begin{remark}\label{claro}
\begin{enumerate}
\item
By \cite[Example 2.6]{Brze1}, a ring $A$ can be viewed as a trivial $A$-coring and the category of comodules of $A$, $\rcomod{A}$ is isomorphic to $\rmod{A}$.
Using this observation, we take $\mathfrak{C} = A$ and $\mathfrak{D} = B$ and we find that Theorem \ref{qF-morphism} reduces to \coref{frobext} where  a functorial
characterization of quasi-Frobenius ring extensions is given. In
this case we have that $A\otimes_AB\cong B$ is $(A,B)$ quasi-finite if and only  $_AB$ is finitely generated and projective. Moreover, $\cohom{B}{B}{-}\simeq -\ot_B\hom{A}{_AB}{A}$.
\item When $A=B$, the corectriction functor $-\ot_AA: \Mm^\cc\to\Mm^\dd$ is quasi-Frobenius if and only if $\cc$ is an $(A,\cc)$ quasi-finite injector and $\cohom{\dd}{\cc}{\dd}\sim \cc$ as $(\dd,\cc)$-bicomodules.
\item
 When $A=B=k$, Theorem \ref{qF-morphism} establishes  that the corectriction functor
$(-)_{\varphi}: \mathcal{M}^C\rightarrow \mathcal{M}^D$ is quasi-Frobenius if and only if $C_D$ is a quasi-finite injector and $\cohom{D}{C}{D}\sim C$ as
$(D,C)$-bicomodules.
\end{enumerate}
\end{remark}
 It is then reasonable to give the following definition.
\begin{definition}\label{defi}
Let $(\varphi, \rho):\mathfrak{C}\rightarrow
\mathfrak{D}$ be a homomorphism of corings such that
$_A\mathfrak{C}$ and $_B\mathfrak{D}$ are flat. It is said to be a
\emph{ right quasi-Frobenius} morphism of corings if
$-\otimes_AB:\mathcal{M}^{\mathfrak{C}}\rightarrow
\mathcal{M}^{\mathfrak{D}}$ is a quasi-Frobenius functor.
\end{definition}

\subsection*{Quasi-Frobenius corings}
Let $\mathfrak{C}$ be an $A$-coring. The forgetful functor
$\functor{U}:\rcomod{\mathfrak{C}}\rightarrow \rmod{A}$ is a left adjoint for the
induction functor $-\otimes_A\mathfrak{C}$ (see \cite[Lemma 3.1]{Brze}).
By classical hom-tensor relations, a right adjoint for the induction functor is given by  $\Hom^\cc(\cc^\cc,-):\Mm^\cc\to\Mm_A$. In light of our previous discussions, it is now a natural question to pose when $\U$ is similar to $\Hom^\cc(\cc^\cc,-)$, i.e.\ when the induction functor is a quasi-Frobenius functor. Comparing this to \coref{frobext}, this leads to the following definition.

\begin{definition}  An $A$-coring $\cc$ is called a \emph{quasi-Frobenius coring} provided the induction functor $-\ot_A\cc:\Mm_A\to \Mm^\cc$ is a quasi-Frobenius functor.
\end{definition}

A characterization of such corings is the following that generalizes \cite[Theorem 4.2]{guo} for left quasi-Frobenius corings.
First, recall that $\*c:={\Hom_A}({_A\cc},{_AA})$ and $\cc^*:=\Hom_A(\cc_A,A_A)$ are rings with unit
$\varepsilon_\cc$ and multiplication for all $f,g\in \*c$ and $f',g'\in\cc^*$ given by the formulas
$$f*g(c)=g(c_{(1)}f(c_{(2)}));\qquad f'*g'(c)=f(g(c_{(1)})c_{(2)}).$$
Furthermore, the maps $i:A\to \*c,\ i(a)(c)=
\varepsilon_{\cc}(c)a$ and $i':A\to\cc^*,\ i'(a)(c)=a\varepsilon_\cc(c)$ are ring morphisms.
There are well-defined functors $\Mm^\cc\to \Mm_{\*c}$ and ${^\cc\Mm}\to {_{\cc^*}\Mm}$ putting
$$m\cdot f=m_{[0]}f(m_{[1]}); \qquad g\cdot n= g(n_{[-1]})n_{[0]};$$
for all $f\in \*c$, $g\in \cc^*$, $m\in M$ and $n\in N$, where $M\in\Mm^\cc$ and $N\in{^\cc\Mm}$.

\begin{theorem}\thlabel{qF-coring} Let $\mathfrak{C}$ be an A-coring
with $_A\mathfrak{C}$ flat. Then the following assertions are equivalent.
\begin{enumerate}[(i)]
\item $\cc$ is a quasi-Frobenius coring;
\item $(\U',\cc\ot_A-,{\Hom^\cc}({^\cc\cc},-))$ is a quasi-Frobenius triple between ${_A\Mm}$ and ${^\cc\Mm}$, where $\U':{^\cc\Mm}\to{_A\Mm}$ is the forgetful functor, i.e.\ $\cc\ot_A-:{_A\Mm}\to{^\cc\Mm}$ is a quasi-Frobenius functor;
\item   $\cc$ is finitely generated projective as a left $A$-module and $\cc \sim \, \*c$ as $(A,\*c)$-bimodules;
\item $\cc$ is finitely generated and projective as a right $A$-module and $\cc\sim \cc^*$ as $(\cc^*,A)$-bimodules;
\item $\cc$ is finitely generated and projective as a left $A$-module and $(-\ot_A\*c,\U,-\ot_A\cc)$ is a quasi-Frobenius triple of functors between $\Mm_A$ and $\Mm^\cc$, i.e.\ $\U:\Mm_A\to\Mm^\cc$ is a quasi-Frobenius functor.
\item $\mathfrak{C}$ is finitely
generated projective as left $A$-module and  $i:A\rightarrow \*c$ is a
quasi-Frobenius extension;
\item $\cc$ is a quasi-Frobenius $(A,\*c)$-bimodule;
\item $\*c$ is a quasi-Frobenius $(\*c,A)$-bimodule
\item `left-right' duals of (v)-(viii), replacing the categories of right (co)modules by their left counterparts and replacing $\*c$ by the ring $\cc^*$.
\end{enumerate}
\end{theorem}
\begin{proof}
$\ul{(i)\Rightarrow(iv)}$. Consider $A$ as a trivial $A$-coring. If $\cc$ is quasi-Frobenius, then we have a pair of quasi-Frobenius functors between the categories $\Mm^A$ and $\Mm^\cc$. Applying Theorem \ref{main}, we find a $(\cc,A)$-bicomodule $X$ that is an $(A,A)$ quasi-finite injector and such that $\Hom^\cc(\cc,-)\cong -\cot_\cc X$. By \cite[Example 4.3]{gomez}, the quasi-finiteness of $X$ implies that $X$ is finitely generated and projective as a left $A$-module.
Moreover, since $X\cong \cc\cot_\cc X\cong \Hom^\cc(\cc,\cc)\cong \cc^*$, we find that $\cc^*$ is finitely generated and projective as a left $A$-module, and therefore $\cc$ is finitely generated as a right $A$-module. Finally, $\U(\cc)\sim \Hom^\cc(\cc,A)$ as right $A$-modules, and by naturality of the functors, we obtain that similarity holds as well as left $\cc$-comodules, hence as left $\cc^*$-modules.\\
$\ul{(iv)\Rightarrow(iii)}$.
Since $\cc$ is finitely generated and projective as a right $A$-module, $\cc^*$ is finitely generated and projective as a left $A$-module. But since $\cc\sim \cc^*$, we find by \leref{divfgp} that $\cc$ is also finitely generated and projective as a left $A$-module. Furthermore, applying the functor $\Hom_A(-,{_AA})$ to $\cc\sim\cc^*$ we obtain that $\*c\sim \Hom_A({_A\cc^*},{_AA})\cong \cc$.\\
$\ul{(iii)\Rightarrow(viii)}$. Since $\cc$ is finitely generated as a left $A$-module, $\*c$ is finitely generated as a right $A$-module. Obviously, $\*c$ is finitely generated as a left $\*c$-module. Moreover $\Hom_A(\*c_A,A_A)\cong \cc$ and $\Hom_{\*c}({_{\*c}\*c},{_{\*c}\*c})\cong \*c$. The implication follows now from \deref{qfbim}\\
\ul{$(viii)\Leftrightarrow(vi)$}. This is a direct application of \coref{frobext}.\\
\ul{$(viii)\Leftrightarrow(v)$}. Consider the following diagram of functors
\begin{equation}\eqlabel{qfcoring}
\xymatrix{
\Mm_A \ar[rr]^-{-\ot_A\*c} \ar@{=}[d] && \Mm_{\*c} \ar[rr]^-{-\ot_{\*c}\*c} \ar[d]^\cong &&
\Mm_A \ar[rr]^-{-\ot_A\cc} \ar@{=}[d] && \Mm_{\*c} \ar[d]^-\cong \ar[rrr]^-{-\ot_{\*c}\Hom_{\*c}(\cc_{\*c},\*c_{\*c})} &&& \Mm_A \ar@{=}[d]\\
\Mm_A \ar[rr]^-{-\ot_A\*c} && \Mm^\cc \ar[rr]^-{\U} &&
\Mm_A \ar[rr]^-{-\ot_A\cc} && \Mm^\cc \ar[rrr]^-{\Hom^\cc(\cc^{\cc},-)} &&& \Mm_A
}
\end{equation}
Since $\cc$ is finitely generated and projective as a left $A$-module, $\Mm^\cc\cong \Mm_{\*c}$, in particular, the functor $-\ot_A\*c:\Mm_A\to\Mm_{\*c}$ makes sense, since $\*c$ is a right $\cc$-comodule. Consider the three first functors in the upper and lower row of the diagram. Clearly, the diagram commutes on this part. On the upper row, the triple is quasi-Frobenius if and only if $\*c$ is a quasi-Frobenius $(\*c,A)$-bimodule, the triple on the lower row will be quasi-Frobenius if and only if condition $(v)$ holds. As the categories are all isomorphic, this shows the equivalence.\\
\ul{$(viii)\Leftrightarrow(vii)$}. This follows from \reref{sequence}.\\
\ul{$(vii)\Rightarrow(i)$}. Consider again diagram \equref{qfcoring}. Condition $(vii)$ means exactly that the last tree functors in the upper row of the diagram are a quasi-Frobenius triple. Furthermore, $\cc$ is finitely generated and projective both as a right $\*c$ module and as a left $A$-module, therefore there are natural isomorphisms
$-\ot_{\*c}\Hom_{\*c}(\cc_{\*c},\*c_{\*c})\simeq \Hom_{\*c}(\cc_{\*c},-)\simeq \Hom^\cc(\cc^\cc,-)$,
this means that the last square of diagram \equref{qfcoring} is commutative (the commutativity of the remaing part was already checked before). Hence the last tree functors of the lower row in the diagram are also a quasi-Frobenius triple, i.e.\ $\cc$ is a quasi-Frobenius coring.\\
Finally, by left-right symmetric arguments, one proves $(ii)\Rightarrow(iii)\Rightarrow(iv)\Rightarrow(ix)\Rightarrow(ii)$.
\end{proof}

\begin{corollary}Let $\varphi:\cc\to\dd$ be a right quasi-Frobenius homorphism of $A$-corings.
If $\dd$ is a quasi-Frobenius $A$-coring, then $\cc$ is also a quasi-Frobenius $A$-coring.
\end{corollary}
\begin{proof}
Consider the following commutative diagram of functors,
\[
\xymatrix{
\Mm^\cc \ar@<-.5ex>[ddrr]_-{U^\cc} \ar@<.5ex>[rrrr]^-{-\ot_AA} &&&& \Mm^\dd \ar@<.5ex>[ddll]^-{U^\dd}
\\  \\
&& \Mm_A
}
\]
where $U^\cc$ and $U^\dd$ denote the forgetful functors.
The statement follows now immediately from Lemma \ref{properties} $(b)$.
\end{proof}

\begin{remark}
A related, more general notion is that of a \emph{quasi-co-Frobenius} coring.
A categorical description of quasi-co-Frobenius corings was initiated recently by Iovanov and the third author in \cite{Iov-Vercr}.
\end{remark}

Recall from \cite{Sw2}, that given a ring extension
  $\rho:R\rightarrow S$, one can view $\mathfrak{C} = S\otimes_RS$ as
  an $S$-coring. This construction is presently known as \emph{Sweedler's coring}
  associated to $\rho$.
The following result generalizes \cite[Theorem 2.7]{Brze1} and \cite[Proposition 4.3]{guo} and can be viewed as the endomorphism ring theorem for quasi-Frobenius extension  in terms of corings.
\begin{proposition} Let  $\mathfrak{C} = S\otimes_RS$ be   the Sweedler's coring associated to a ring extension $\rho:R\rightarrow S.$ If $S$ is a quasi-Frobenius extension of $R$, then $\mathfrak{C}$ is a quasi-Frobenius $S$-coring.
\end{proposition}
\begin{proof}
If $\rho$ is a quasi-Frobenius ring extension, then we know by \coref{frobext} that $\rho_*$ is a quasi-Frobenius functor, $_RS$ is finitely generated projective and $_SS_R\sim (_RS)^*$ as $(S,R)$-bimodules.
By Lemma \ref{properties}, the functor $-\otimes_RS$  preserves finitely generated and projective modules. Hence, $(S\otimes_RS)_S$ is finitely generated and projective as a right $S$-module. Applying $-\otimes_RS$ to $_SS_R\sim (_RS)^*$ we obtain
$$_R(S\otimes_RS)_S\sim (_RS)^*\otimes_RS\cong \End{R}{}{S}.$$
Now from \cite[Proposition 2.1]{kaouti-pepe}, $\End{R}{}{S}\cong\overline{T}$, where $\overline{T}$ is the opposite algebra of $((S\otimes_RS)_S)^*$. Therefore $S\otimes_RS\sim \overline{T},$ and $S\otimes_RS$ is a quasi-Frobenius $S$-coring by \thref{qF-coring}.
\end{proof}

\bigskip

\bibliographystyle{amsplain}

\end{document}